\documentclass[12pt,draft]{amsart}
\usepackage{amssymb,latexsym, amsmath, amscd, array, graphicx}

\newtheorem{thm}{Theorem}[section]
\newtheorem{theorem}[thm]{Theorem}
\newtheorem{lem}[thm]{Lemma}
\newtheorem{lemma}[thm]{Lemma}

\newtheorem{prop}[thm]{Proposition}

\theoremstyle{definition}

\theoremstyle{definition}

\newtheorem{definition}[thm]{Definition}

\newtheorem{remark[thm]}{Remark}

\newtheorem{cor}[thm]{Corollary}

\theoremstyle{remark}

\DeclareMathOperator{\cat}{{\mbox{\rm cat$_{\rm LS}$}}}

\def\Ord{\protect\operatorname{Ord}}

\def\scr{\mathcal}

\def\Z{{\mathbb Z}}

\def\R{{\mathbb R}}

\long\def\forget#1\forgotten{} %

\numberwithin{equation}{section}

\begin{document}

\title[Lusternik-Schnirelmann category]{
The Lusternik-Schniremann-category and the fundamental group}
\author[A.~Dranishnikov]{Alexander N. Dranishnikov}


\address{Department of Mathematics, University
of Florida, 358 Little Hall, Gainesville, FL 32611-8105, USA}
\curraddr{}
\email{dranish@math.ufl.edu}
%\thanks{Supported by NSF, grant DMS-0904278}

%    \subjclass is required.

\date{}

\dedicatory{}

%    "Communicated by" -- provide editor's name; required.
\commby{Daniel Ruberman}

%    Abstract is required.
\begin{abstract}
We prove that
$$
\cat X\le cd(\pi_1(X))+\bigg\lceil\frac{\dim X-1}{2}\bigg\rceil$$
for every CW complex $X$ where $cd(\pi_1(X))$ denotes the cohomological
dimension of the fundamental group of $X$.
\end{abstract}

\maketitle

\section{Introduction} The {\em Lusternik-Schnirelmann category}  
(briefly LS-category) $\cat X$ of a
topological space $X$ is the minimal number $n$ such that there is an open cover
$\{U_0,\dots, U_n\}$ of $X$ by $n+1$ contractible in $X$ sets. We
note that sets $U_i$ are not necessarily contractible. The
Lusternik-Schnirelmann category proved to be useful in different
areas of mathematics. In particular, the classical theorem of
Lusternik and Schnirelmann \cite{CLOT} proven in the 30s states that
$\cat M$ gives a low bound for critical points on $M$ of any smooth
not necessarily Morse function. For nice spaces as CW complexes it
is an easy observation that $\cat X\le\dim X$. In the 50s Whitehead
proved that for simply connected CW complexes $\cat X\le\dim X/2$
\cite{CLOT}. In a presence of the fundamental group as small as
$\Z_2$ the Lusternik-Schnirelmann category can be equal the
dimension. An example is $\R P^n$.

Nevertheless, Yu.~Rudyak conjectured that in the case of free
fundamental group there should be Whitehead's type inequality at
least for closed manifolds. There were partial results towards
Rudyak's conjecture \cite{DKR},\cite{St} until it was settled in
\cite{Dr}. Also it was shown in \cite{Dr} that Whitehead's type
estimate holds for complexes with the fundamental group of
cohomological dimension $\le 2$. We recall that free groups (and
only them  \cite{Stal},\cite{Swan})
have cohomological dimension one. In this paper we prove Whitehead's
type inequality for complexes with fundamental groups having finite
cohomological dimension.

We conclude the introductory part by definitionss and 
statements from \cite{Dr} which
are used in this paper.
Let $\scr U=\{U_{\alpha}\}_{\alpha\in A}$ be a family
of sets in a topological space $X$.
Formally, it is a function $U: A\to 2^X\setminus\{\varnothing\}$ form
the index set to the set of nonempty subsets of $X$.
The sets $U_{\alpha}$ in the family $\scr U$ will be called
{\em elements of} $\scr U$.
The {\em multiplicity} of $\scr U$ (or the {\em order}) at a point $x\in
X$, denoted $\Ord_x\scr U$, is the number of elements of $\scr U$
that contain $x$. The {\em multiplicity} of $\scr U$ is defined as
$\Ord\scr U=\sup_{x\in X}\Ord_x\scr U$. A family $\scr U$
is a cover of $X$ if $\Ord_x\scr U\ne 0$ for all $x$.
A cover $\scr U$ is a {\em
refinement} of another cover $\scr C$ ($\scr U$ {\em refines} $\scr
C$) if for every $U\in\scr U$ there exists $C\in\scr C$ such that
$U\subset C$. We recall that the {\em covering dimension} of a
topological space $X$ does not exceed $n$, $\dim X\le n$, if for
every open cover $\scr C$ of $X$ there is an open refinement $\scr
U$ with $\Ord\scr U\le n+1$.

\begin{definition} A family $\scr U$ of subsets of $X$ is called a {\em
$k$-cover}, $k\in N$ if every subfamily of $k$ elements forms a
cover of $X$.
\end{definition}
The following is obvious (see \cite{Dr}).
\begin{prop}\label{n-cover}
A family $\scr U$ that consists of $m$ subsets of $X$ is an
$(n+1)$-cover of $X$ if and only if $\Ord_x\scr U\ge m-n$ for all
$x\in X$.
\end{prop}

The following theorem can be found in \cite{Os}.
\begin{theorem}[Kolmogorov-Ostrand]\label{Ost}
A metric space $X$ is of dimension $\leq n$ if and only if for each
open cover $\scr C$ of $X$ and each integer $m\geq n+1$, there exist
$m$ disjoint families of open sets $\scr U_1,\cdots,\scr U_m$ such
that their union $\cup\scr U_i$ is an  $(n+1)$-cover of $X$ and it
refines $\scr C$.
\end{theorem}

Let $f:X\to Y$ be a map and let $X'\subset X$. A set $U\subset X$ is
{\em fiber-wise contractible to $X'$} if there is a homotopy
$H:U\times[0,1]\to X$ such that $H(x,0)=x$, 
$H(U\times\{1\})\subset X'$, and $f(H(x,t))=f(x)$ for all $x\in U$.

We refer to \cite{Dr} for the proof of the following
\begin{theorem}\label{criterion}
Let $\scr U=\{U_0,\dots,U_n\}$ be an open cover of a normal
topological space $X$. Then  for any $m=n,n+1,\dots,\infty$ there is
an  open $(n+1)$-cover of $X$, $\scr U_m=\{U_0,\dots,U_m\}$
extending $\scr U$ such that  for $k>n$,
 $U_k=\cup_{i=0}^nV_i$ is a disjoint union with
$V_i\subset U_i$.
\end{theorem}

\begin{cor}\label{collcriterion}
Let $f:X\to Y$ be a continuous map of a normal topological space and
let $\scr U=\{U_0,\dots,U_n\}$ be an open cover of  $X$ by sets
fiber-wise contractible to $X'$. Then for any $m=n,n+1,\dots,\infty$
there is an open $(n+1)$-cover of $X$, $\scr U_m=\{U_0,\dots,U_m\}$
by sets fiber-wise contractible to $X'$.
\end{cor}

\section{Generalization of Ganea's fibrations}

Let $A\subset Z$ be a closed subset. By $P_AZ$ we denote the space
of paths issued from $A$, i.e. the space of continuous maps
$\phi:[0,1]\to Z$ with $\phi(0)\in A$. The we define a map
$p_A:P_AZ\to Z$ by the formula $p(\phi)=\phi(1)$. Clearly, $p_A$ is
a Hurewicz fibration. Let $F$ be its fiber.
\begin{prop}\label{fiber}
There is a Hurewicz fibration $\pi:F\to A$ with the fiber $\Omega Z$,
the loop space on $Z$.
\end{prop}
\begin{proof}
The map $q':P_AZ\to A\times Z$ that sends a path to the end points 
is a Hurewicz fibration as a pull-back of the Hurewicz fibration
$q:Z^{[0,1]}\to Z\times Z$ \cite{Sp}. The fiber of $q$ is the loop space
$\Omega Z$. Since $p_A=pr_2\circ q'$, the fiber 
$F=p_A^{-1}(x)=(q')^{-1}pr_2^{-1}(x)=q^{-1}(A)$ is the total space of 
a Hurewicz fibration  $q$ over $A$ with the fiber $\Omega Z$.
\end{proof}

We define the $k$-th {\em generalized Ganea's fibration}
$p_k:E_k(Z,A)\to Z$ over a path connected space $Z$ with a fixed
closed subset $A$ as the fiber-wise join product of $k+1$ copies of
the fibrations $p_A:PZ\to Z$.  Since $p_A$ is a Hurewicz
fibration and the fiber-wise join of Hurewicz fibrations is a
Hurewicz fibration, so are all $p_k$ \cite{Sv}. 
Note that the fiber of $p_k$ is
the join product $\ast^{k+1}F$ of $k+1$ copies of $F$ (see
\cite{CLOT} for more details). Also we note that
for $A=\{z_0\}$ the fibration $p_k$ is the standard Ganea fibration.
The following is a generalization of
the Ganea-\v Svarc theorem.
\begin{thm}\label{generalization} Let $A\subset X$ be a subcomplex
contractible in $X$. Then $\cat (X)\le k$ if and only if the
generalized Ganea fibration $$p_k:E_k(Z,A)\to Z$$ admits a section.
\end{thm}
\begin{proof}
When $A$ is a point this statements turns into the classical
Ganea-\v Svarc theorem (\cite{CLOT}, \cite{Sv}). Since for $z_0\in
A$, the fibration $p_k:E_k(Z,z_0)\to Z$ is contained in
$p_k:E_k(Z,A)\to Z$, the classical Ganea-\v Svarc theorem implies
the only if direction.

The barycentric coordinates of a section to $p_k$ define an open
cover $U_0,\dots U_k$ of $U_i$ with each $U_i$ contractible to $A$.
Since $A$ is contractible in $Z$, all sets $U_i$ are contractible in
$Z$.
\end{proof}
We call a map $f:X\to Y$ a {\em stratified locally trivial bundle}
(with two strata) with fiber $(Z,A)$ if there $X'\subset X$, such that 
$(f^{-1}(y),g^{-1}(y))\cong
(Z,A)$ for all $y\in Y$, where $g=f|_{X'}$, and there is an open 
cover $\scr U=\{U\}$ of $Y$ such that
$(f^{-1}(U),g^{-1}(U))$ is homeomorphic as a pair to $(Z\times U,
A\times U)$ by means of a fiber preserving homeomorphism. Such a 
bundle is called a {\em trivial startified bundle} if one cant take $\sc U$
consisting of one element $U=Y$.

Now let $f:X\to Y$ be a stratified locally trivial bundle with a
subbundle $g:X'\to Y$ and a fiber $(Z,A)$. We define a space
$$
E_0=\{\phi\in C(I,X)\mid  f\phi(I)=f\phi(0), \ \phi(0)\in
g^{-1}(f\phi(0))\}
$$
to be the space of all paths $\phi$ in $f^{-1}(y)$ for all $y\in Y$
with the initial point in $g^{-1}(y)$. The topology in $E_0$ is
inherited from $C(I,X)$. We define a map $\xi_0:E_0\to X$ by the
formula $\xi_0(\phi)=\phi(1)$. Then $\xi_k:E_k\to X$ is defined as
the fiber-wise join of $k+1$ copies of $\xi_0$. Formally, we define
inductively $E_{k}$ as a subspace of the join $E_0\ast E_{k-1}$:
$$
E_{k}=\cup\{\phi\ast\psi\in E_0\ast E_{k-1}\mid
\xi_0(\phi)=\xi_{k-1}(\psi)\}
$$
which is the union of all intervals $[\phi,\psi]=\phi\ast\psi$ with
the endpoints $\phi\in E_0$ and $\psi\in E_{k-1}$ such that
$\xi_0(\phi)=\xi_{k-1}(\psi)$. There is a natural projection
$\xi_{k}:E_{k}\to X$ that takes all points of each interval
$[\phi,\psi]$ to $\phi(0)$.

Note that when $f:X=Z\times Y\to Y$ is a trivial stratified bundle
with the subbundle $g:A\times Y\to Y$, $A\subset Z$, then
$E_k=E_k(Z,A)\times Y$ and $\xi_k=p_k\times 1_Y$ where
$p_k:(E_k,A)\to Z$ is the generalized Ganea fibration.

\begin{lem}\label{Generalized Ganea}
Let $f:X\to Y$ be a stratified locally trivial bundle between
paracompact spaces with a fiber $(Z,A)$ in which $A$ is contractible
in $Z$. Then

i. For each $k$ the map $\xi_k:E_k\to X$ is a Hurewicz fibration.

ii. The fiber of $\xi_k$ is  the join of $k+1$ copies of the fiber $F$
of $p_A:P_AZ\to Z$.

iii. If the projection $\xi_k$ has a section, then $X$ has an
open cover $\scr U=\{U_0,\dots,U_k\}$ by sets each of which admits a
fiber-wise deformation into $X'$ where $g:X'\to Y$ is the subbundle.
\end{lem}
\begin{proof}
i. First, we note that this statement holds true for trivial
stratified bundles. By the assumption there is a cover $\scr U$ of
$Y$ such that $f|_{f^{-1}{U}}:f^{-1}(U)\to U$ is a trivial stratified
bundle and hence $\xi_k$ is a Hurewicz fibration over $f^{-1}(U)$
for all $U\in\scr U$. Then we apply Dold's theorem \cite{Do} to
conclude that that $\xi_k$ is a Hurewicz fibration over $X$.

ii. We note that $\xi_k$ over $f^{-1}(y)$ coincides with
the generalized Ganea fibration $p_k$ for $(Z,A)$. Therefore, the
fiber of $\xi_k$ coincides with the fiber of $p_k$. Then we apply 
Proposition~\ref{fiber}

iii.  Suppose $\xi_k$ has a section
$\sigma:X\to E_k$. For each $x\in X$ the element $\sigma(x)$ of
$\ast^{k+1}\Omega F$ can be presented as the $(k+1)$-tuple
$$\sigma(x)=((\phi_0,t_0),\dots,(\phi_k,t_k))\mid \sum t_i=1, t_i\ge
0).$$ We use the notation $\sigma(x)_i=t_i$. Clearly, $\sigma(x)_i$
is a continuous function.

A section $\sigma:X\to E_k$ defines a cover $\scr U=\{U_0,\dots,
U_k\}$ of $X$ as follows:
$$
U_i=\{x\in X\mid\sigma(x)_i>0\}.
$$
By the construction of $U_i$ for $i\le n$ for every $x\in U_i$ there
is a canonical path connecting $x$ with $g^{-1}f(x)$. These paths 
define a fiber-wise deformation of $U_i$ into $g^{-1}f(U_i)\subset X'$.
\end{proof}

\section{The main result}

We recall that the {\em homotopical dimension} of a space $X$,
$hd(X)$, is the minimal dimension of a $CW$-complex homotopy
equivalent to $X$ \cite{CLOT}.
\begin{prop}\label{section}
Let $p:E\to X$ be a fibration with $(n-1)$-connected fiber where
$n=hd(X)$. Then $p$ admits a section.
\end{prop}
\begin{proof}
Let $h:Y\to X$ be a homotopy equivalence with the homotopy inverse
$g:X\to Y$ where $Y$ is a CW complex of dimension $n$. Since the
fiber of $p$ is $(n-1)$-connected, the map $h$ admits a lift
$h':Y\to E$. Let $H$ be a homotopy connecting $h\circ g$ with $1_X$.
By the homotopy lifting property there is a lift $H':X\times I\to E$
of $H$ with $H|_{X\times\{0\}}=h'\circ g$. Then the restriction
$H|_{X\times\{1\}}$ is a section.
\end{proof}

We recall that $\lceil x\rceil$ denotes the smallest integer $n$
such that $x\le n$.
\begin{lemma}\label{r-connectedhur}
Suppose that a stratified locally trivial bundle $f:X\to Y$ with a
fiber $(F,A)$ is such that $F$ is  $r$-connected, $A$ is
$(r-1)$-connected,  $A$ is contractible in $F$, and $Y$ is locally
contractible. Then
$$\cat X\le \dim Y+
\bigg\lceil\frac{hd(X)-r}{r+1}\bigg\rceil.$$
\end{lemma}
\begin{proof}
Let $\cat Y=m$ and  $hd(X)=n$.

We note that the fiber $K$ of $p_k:E_k(F,A)\to F$ is the join
product $\ast^{k+1}K_0$ of $k+1$ copies of the fiber $K_0$ of the
map $p_A:PF\to F$. By Proposition~\ref{fiber}, $K_0$ admits a fibration
$\phi:K_0\to A$ with fibers homotopy equivalent to the loop space
$\Omega F$. Since the base $A$ and the fibers are $(r-1)$-connected,
$K_0$ is $(r-1)$-connected. Thus, $K$ is $(k+(k+1)r-1)$-connected.
By Proposition \ref{section} there is a section $\sigma:X\to E_k$
whenever $k(r+1)+r\ge n$. The smallest such $k$ is equal to
$\lceil\frac{n-r}{r+1}\rceil$.

By Lemma~\ref{Generalized Ganea} a section $\sigma:X\to E_k$ defines
a cover $\scr U=\{U_0,\dots, U_k\}$ by the sets fiber-wise
contractible to $X'$ where $X'\subset X$ is the first stratum. Let
$\scr U_{m+k}=\{U_0,\dots, U_{m+k}\}$ be an extension of $\scr U$ to
a $(k+1)$-cover of $X$ from Corollary~\ref{collcriterion}.

Let $\scr V=\{V_0,\dots, V_{m+k}\}$ be an open $(m+1)$-cover of $Y$
such that for every $i$, $V_i$ is contractible to a point in $V_i'$
and $f$ is trivial stratified bundle over $V_i'$. Such $\scr V$
exists in view of Theorem~\ref{Ost}. We show that the sets
$W_i=f^{-1}(V_i)\cap U_i$ are contractible in $X$ for all 
$i\in\{0,1,\dots,m+k\}$. By
Corollary~\ref{collcriterion} $U_i$ is fiber-wise contractible into
$X'$ for $i\le m+k$. Hence we can contract $f^{-1}(V_i)\cap U_i$ to
$f^{-1}(V_i)\cap X'\cong V_i\times A$ in $X$. Then we apply a
contractions to a point of $V_i$ in $V_i'$ and $A$ in $F$ to obtain
a contraction to a point of $f^{-1}(V_i)\cap X'\cong V_i\times A$ in
$f^{-1}(V_i')=V_i'\times F$.

Next we show that $\{W_i\}_{i=0}^{m+k}$ is a cover of $X$. Since $\scr V$
is an $(m+1)$-cover, by Proposition \ref{n-cover} every $y\in Y$ is
covered by at least $k+1$ elements $V_{i_0},\dots, V_{i_k}$ of $\scr
V$. Since $\scr U_{m+k}$ is a $(k+1)$-cover, $U_{i_0},\dots,
U_{i_k}$ is a cover of $X$. Hence $W_{i_0},\dots, W_{i_k}$ covers
$f^{-1}(y)$.
\end{proof}

\begin{theorem}
For every CW complex $X$ with  the following inequality holds true:
$$
\cat X\le cd(\pi_1(X)) +\left\lceil\frac{hd(X)-1}{2}\right\rceil.$$
\end{theorem}
\begin{proof}
Let $\pi=\pi_1(X)$ and let $\tilde X$ denote the universal cover of $X$.
We consider Borel's construction
$$
\begin{CD}
\tilde X @<<< \tilde X\times E\pi @>>> E\pi\\
@VVV @VVV @VVV\\
X @<g<< \tilde X\times_{\pi}E\pi @>f>> B\pi.\\
\end{CD}
$$
The 1-skelwton $X^{(1)}$ of $X$ defines a $\pi$-equivariant 
stratification $\tilde X^{(1)}\subset \tilde X$ of the universal cover.
This startification allows us to
treat $f$ as a stratified locally trivial bundle with the fiber
$(\tilde X,\tilde X^{(1)})$. We note that all condition of
Lemma~\ref{r-connectedhur} are satisfied for $r=1$. Therefore,
$$
\cat X\le \dim B\pi
+\left\lceil\frac{hd(\tilde X\times_{\pi}E\pi)-1}{2}\right\rceil.$$ Since
$g$ is a fibration with the homotopy trivial fiber, the space
$\tilde X\times_{\pi}E\pi$ is homotopy equivalent to $X$. Thus,
$hd(\tilde X\times_{\pi}E\pi)=hd(X)$. In view of the results of
Eilenberg and Ganea \cite{EG},\cite{Br} we may assume that $\dim
B\pi=cd(\pi)$ if $cd(\pi)>2$. The case when $cd(\pi)\le 2$ is
treated in \cite{Dr}.
\end{proof}




\begin{thebibliography}{[CLOT]}
\bibliographystyle{amsplain}

\bibitem[CLOT]{CLOT} Cornea, O.; Lupton, G.; Oprea, J.; Tanr\'e, D.:
Lusternik-Schnirelmann category.  {\em Mathematical Surveys and
Monographs}, \textbf{103}.  American Mathematical Society, Providence,
RI, 2003.

\bibitem[Br]{Br} Brown, K.: Cohomology of groups, Springer, 1982.

\bibitem[Do]{Do} Dold, A.: Partition of unity in the theory of fibrations.
{\em Ann. of Math.}, \textbf{78} (1963), 223-255.


\bibitem[DKR]{DKR} Dranishnikov, A; Katz, M.; Rudyak, Y.:
Small values of Lusternik-Schnirelmann
category and systolic categories for manifolds.  {Preprint},  (2007)
\texttt{arXiv:0706.1625math.AT} (to appear in Geometry and Topology).

\bibitem[Dr]{Dr} Dranishnikov, A.:On the Lusternik-Schnirelmann 
category of spaces with 2-dimensional 
fundamental group.  Proc. Amer. Math. Soc.  137  (2009),  no. 4, 1489--1497. 

\bibitem[EG]{EG} Eilenberg, Samuel; Ganea, Tudor: On the Lusternik-Schnirelmann
category of abstract groups.
{\em Ann. of Math.} (2)  \textbf{65}  (1957), 517--518.



\bibitem[Os]{Os} Ostrand, Ph. : Dimension of metric spaces and Hilbert's
problem $13$, Bull. Amer. Math. Soc. 71 1965, 619-622.

\bibitem[Sp]{Sp} Spanier, E.: Algebraic Topology, Mebman - Hell, 1966.


\bibitem[Sta]{Stal} Stallings, J.: Groups of dimension 1 are locally
free. \emph{Bull. Amer. Math. Soc.} \textbf{74} (1968), 361--364.


\bibitem[St]{St} Strom, J.: Lusternik-Schnirelmann category of
spaces with free fundamental group.   {\em AGT}, \textbf{7}  2007, 1805-1808.

\bibitem[Sv]{Sv} \v Svarc, A.: The genus of a fibered space. {\em Trudy
Moskov. Mat. Ob\v s\v c} \textbf {10, 11} (1961 and 1962), 217--272, 99--126,
(in {\em Amer. Math. Soc. Transl.} Series 2, vol \textbf{55} (1966)).

\bibitem[Swan]{Swan} Swan, R.: Groups of cohomological dimension
one.  {\em J. Algebra\/} \textbf{12} (1969), 585--610.


\end{thebibliography}
\end{document}